# The Complexity of Broadcasting in Bounded-Degree Networks

Michael J. Dinneen
Computer Research and Applications
Los Alamos National Laboratory
Los Alamos, New Mexico 87545, U.S.A.

## 1  Introduction

Because of cost and physical limitations there is an interest in broadcast networks with a fixed maximum number of communication links for each node, [BHLP,LP,DFF]. In this note we show that given an arbitrary bounded-degree network, determining the minimum broadcast time for any originating node is NP-complete. For the NP-completeness proof of the general case of unbounded-degree networks, it is straightforward to reduce the instances of the three dimensional matching problem (3DM) to a corresponding minimum broadcast time problem (MBT), [GJ].

Broadcasting concerns the dissemination of a message originating at one node of a network to all other nodes, [FHMP]. This task is accomplished by placing a series of calls over the communication lines of the network between neighboring nodes, where each call requires a unit of time and a call can involve only two nodes. A node can participate in only one call per unit of time.

In the next section we introduce terminology and present some initial observations. In the second section we show that for bounded-degree networks determining the minimum broadcast time from an originating node remains NP-complete. In the last section we propose a problem regarding bounded-degree Cayley broadcast graphs.

## 2  Preliminaries

Graphs in this paper are simple and undirected. We use the terms network and graph interchangeably, and similarly node and vertex. Let $G = (V, E)$ be a connected graph and let $u$ be a vertex of $G$. The **broadcast time of vertex** $u$, $b(u)$, is the minimum number of time units required to complete broadcasting of a message originating at vertex $u$. The **broadcast time of** $G$ is the maximum broadcast time of any vertex $u$ in $G$, $b(G) = \max\{b(u) \mid u \in V\}$. The decision problem that we address is the following where the integer $C$ is some bound on the vertex degree.

**<u>Bounded-Degree Minimum Broadcast Time Problem</u>**

> Instance: Graph $G = (V, E)$ with maximum degree $C$, vertex $v \in V$, and a positive integer $k$.
> Question: Is the broadcast time $b(v)$ less than or equal to $k$? Or explicitly, is there a sequence $V_0 = \{v\}, E_1, V_1, E_2, \ldots, E_k, V_k$ such that each $V_i \subseteq V$, each $E_i \subseteq E$, $V_k = V$, and, for $1 \leq i \leq k$, (1) each edge in $E_i$ has exactly one endpoint in $V_{i-1}$, (2) no two edges in $E_i$ share a common endpoint, and (3) $V_i = V_{i-1} \cup \{v \mid (u, v) \in E_i\}$ ?

One graph that is a building block for our NP-complete reduction is the following (see Figure 1). The graph $A_n = (V_n, E_n)$ where $V_n = \{(i, j) \mid 1 \leq i \leq n$ and $1 \leq j \leq 2(n - i) + 1\}$ and $E_n = \{((i, 1), (i + 1, 1)) \mid 1 \leq i \leq n - 1\} \cup \{((i, j), (i, j + 1)) \mid 1 \leq i \leq n - 1$ and $1 \leq j \leq 2(n - i)\}$. For convienence, vertex $(1, 1)$ is labeled $r$ and vertices $(i, 2n - 2i + 1)$ are labeled $v_i$ in this graph. The graph $A_n$ is a binary tree (not complete) with $\sum_{k=1}^{n} 2k - 1 = n^2$ vertices.

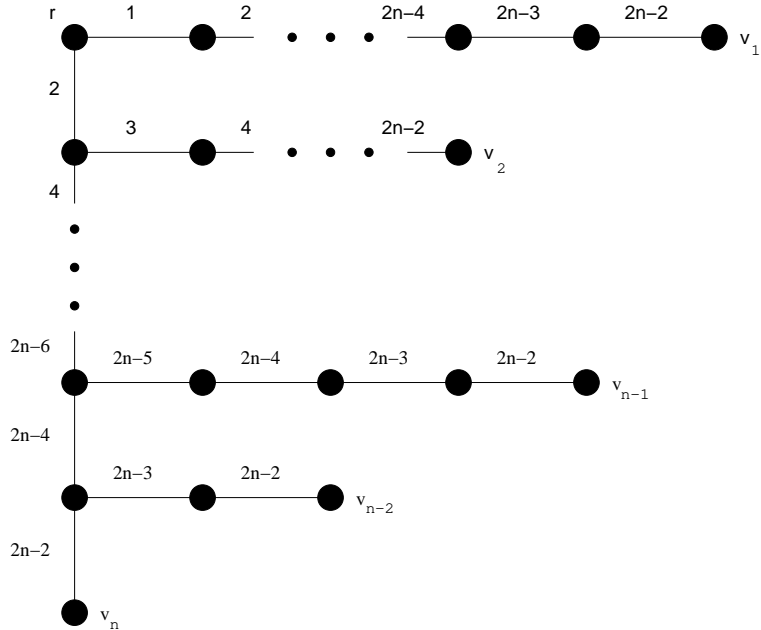

Figure 1: The graph $A_n$.

**Lemma 1**: The broadcast time of vertex $r$ in graph $A_n$ is $2n - 2$.

To illustrate Lemma 1, the edges in Figure 1 are labeled to show the unique routing scheme necessary to achieve the minimum broadcast time.

The next two lemmas assume that $A_n$ ($n > 1$) is a vertex-induced subgraph of some graph $G$, only vertices $\{r, v_1, v_2, \ldots, v_n\}$ have neighbors outside $A_n$, and $r$ receives a message at time $t$.

**Lemma 2**: Every vertex $v_i$ of $A_n$ can simultaneous relay the broadcast message out of $A_n$ at time $t + 2n - 1$.

**Lemma 3**: If some vertex $v_i$ of $A_n$ relays a message out of $A_n$ at time less than $t + 2n - 1$ then the broadcast time locally in $A_n$ is not minimal.

## 3 Main result

**Theorem**: Bounded-Degree Minimum Broadcast Time (BDMBT) is NP-complete.

**Proof.** It is easy to see that BDMBT is in NP since given a routing scheme for a graph $G = (V, E)$ consisting of the sets $V_0 = \{v\}, E_1, V_1, E_2, \ldots, E_k, V_k$ one can determine that each $V_i \subseteq V$, each $E_i \subseteq E$, $V_k = V$, and, for $1 \leq i \leq k$, (1) each edge in $E_i$ has exactly one endpoint in $V_{i-1}$, (2) no two edges in $E_i$ share a common endpoint, and (3) $V_i = V_{i-1} \cup \{v \mid (u, v) \in E_i\}$ in linear time.

We will prove BDMBT is NP-complete by restriction to 3SAT (see [GJ]). That is, given any instance of 3SAT we will construct an instance of BDMBT such that 3SAT has a truth assignment if and only if a specific vertex of BDMBT can broadcast within some time limit.

Given an instance $E$ of 3SAT with $n$ variables and $m$ clauses we build the following graph as depicted in Figure 2. This graph $G$ contains one copy of $A_n$, $2n$ copies of $A_m$, and $m$ additional vertices. For convienience, the root vertices of the $A_m$ graphs are labeled $T_i$ and $F_i$ for each of the $n$ variables $v_i, 1 \leq i \leq n$. Vertex $v_i$ of $A_n$ is attached to $T_i$ and $F_i$ for all $1 \leq i \leq n$. For each of the three literals in a clause $j$, vertex $j$ is attach to vertex $v_j$ of $T_i$'s or $F_i$'s $A_n$ graph depending on whether the literal is positive or negative and has $v_i$ as a variable.

We claim that the graph $G$ has broadcast time $2n + 2m - 2$ from vertex $r$ if and only if the 3SAT instance $E$ has a truth assignment that evaluates to true.

First assume that $E$ is satisfiable. In $2n - 2$ time steps the vertices $v_i$ of $A_n$ for all $1 \leq i \leq n$ can all broadcast to a neighbor outside of $A_n$. At the next broadcast time vertex $v_i$ relays its first message to $T_i$ or $F_i$ depending on whether this variable is set true or false in $E$. After $2m - 2$ more broadcasts, all of the vertices in one half of the $A_m$ graphs will have seen the message, and the other half of the $A_m$ graphs will be within one broadcast of being in the same state. At time $2n + 2m - 2$, every vertex $j$ for $1 \leq j \leq m$ can receive the message from any one of the $A_m$ graphs that corresponds to setting some literal true in clause $j$. Thus, this routing scheme can broadcast in $2n + 2m - 2$ time steps from vertex $r$ if the expression $E$ is satisfiable.

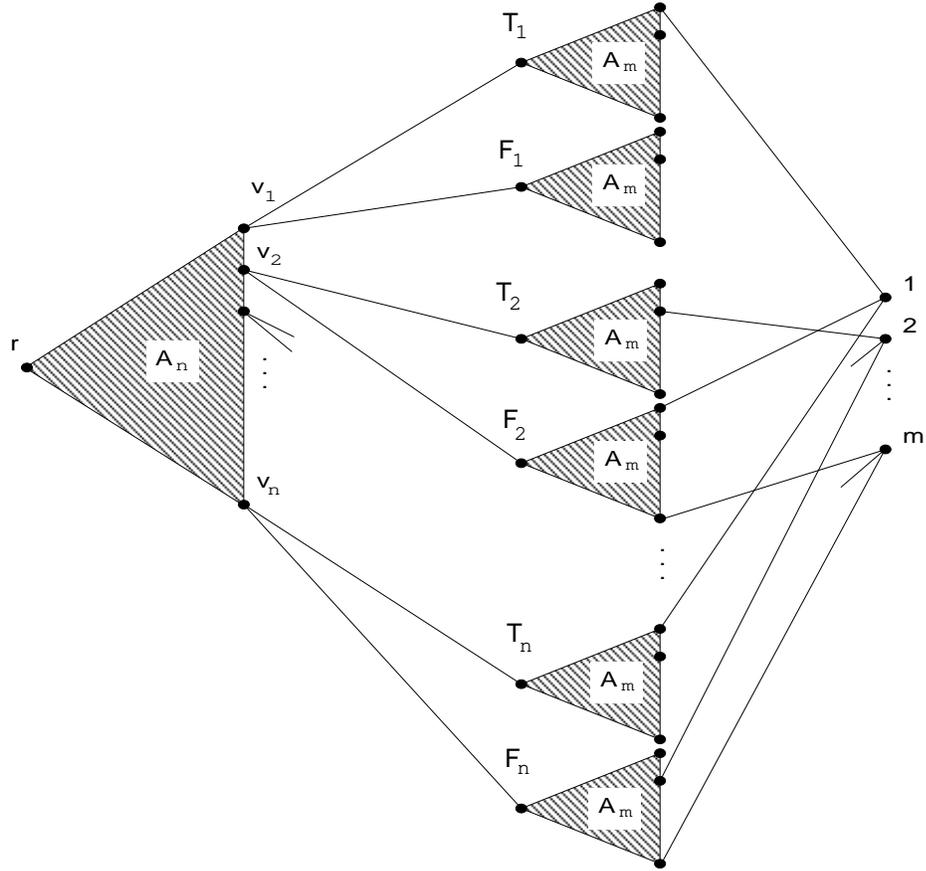

Figure 2: A BDMBT graph for a $n$-variable $m$-clause 3SAT instance.

Now assume that vertex $r$ of the graph $G$ can broadcast in $2n + 2m - 2$ steps and the expression $E$ is not satisfiable. Clearly the above routing scheme was not used since some clause-vertex $c$ will always be adjacent to neighbors that all receive the broadcast message at time $2n + 2m - 2$. Investigating the possibility that the unique minimal routine scheme was not used for $A_n$ reveils that one of the $T_i$'s or $F_i$'s $A_m$ graphs must contain a vertex with shortest distance greater than $2m-1$ from all vertices that could have received the broadcast message at time $2n-2$. Henceforth, routing in $A_n$ must be minimal. Now assume that a minimal routing scheme was used for $A_n$ and an arbitrary truth assignment has been used so that clause $c$ is false (i.e., using the appropriate edges $(v_i, T_i)$ or $(v_i, F_i)$ for each $v_i$). If one of the neighbors of vertex $c$ received the message at time $2n + 2m - 2$ or earlier, then from Lemma 3 this $A_m$ graph is not minimal and one of its $v_i$ vertices will receive the message after time $2n + 2m - 2$. (Notice that the root of this $A_m$ receives the message at time $2n$.) These contradictions show that the minimal broadcast time must be greater than $2n + 2m - 2$. □

## 4 Final Remarks

We have just shown that for any originating vertex in a graph with maximum degree three determining the minimum time needed to broadcast a message is NP-complete. For any maximum degree 2 graph $G$ (a cycle or a path) the answer is trivial (broadcast time of $\lfloor \frac{|G|}{2} \rfloor$ or $|G| - 1$, respectively). Broadcasting in maximum degree 3 directed graphs is also NP-complete by simply directing the edges away from vertex $r$ in the above proof.

The broadcast time $b(G)$ of a given graph $G = (V, E)$ is only slightly harder than the above problem since by definition it is simply the maximum broadcast time $b(u)$ for all $u \in V$. For all vertex-symmetric graphs, in particular Cayley graphs, the broadcast time of each vertex is the same.

Many of the largest graphs with fixed degree and broadcast time are known to be Cayley graphs. Since a general-purpose efficient algorithm does not seem to exist, an open question arrises as of whether the group structure can be exploited in this class of graphs to give the minimum broadcast time directly.

Two examples of Cayley graphs with known simple routing algorithms that yield their minimum broadcast time are the hypercubes and dihedral broadcast graphs. In light of the complexity of BDMBT, heuristics were used to to find short routing schemes in some of the other Cayley broadcast graphs. (See [DFF].)

## References

[BHLP] J-C. Bermond, P. Hell, A. L. Liestman and J. G. Peters, "Broadcasting in Bounded Degree Graphs," Technical Report CMPT 88-5, School of Computing Science, Simon Fraser University, B. C., Canada, 1988.